\documentclass[a4paper, 11pt]{amsart}
\usepackage{graphicx}
\usepackage{amsfonts}
\usepackage{amsmath}
\usepackage{amssymb}
\usepackage{hyperref}

\begin{document}

\title{Quantum computing with Bianchi groups}

\author{Michel Planat$\dag$, Raymond Aschheim$\ddag$,\\ Marcelo M. Amaral$\ddag$ and Klee Irwin$\ddag$}

\address{$\dag$ Universit\'e de Bourgogne/Franche-Comt\'e, Institut FEMTO-ST CNRS UMR 6174, 15 B Avenue des Montboucons, F-25044 Besan\c con, France.}
\email{michel.planat@femto-st.fr}

\address{$\ddag$ Quantum Gravity Research, Los Angeles, CA 90290, USA}
\email{raymond@QuantumGravityResearch.org}
\email{Klee@quantumgravityresearch.org}
\email{Marcelo@quantumgravityresearch.org}

\begin{abstract}

It has been shown that non-stabilizer eigenstates of permutation gates are appropriate for allowing $d$-dimensional universal quantum computing (uqc) based on minimal informationally complete POVMs. The relevant quantum gates may be built from subgroups of finite index of the modular group $\Gamma=PSL(2,\mathbb{Z})$ [M. Planat, Entropy 20, 16 (2018)] or more generally from subgroups of fundamental groups of $3$-manifolds [M. Planat, R. Aschheim, M.~M. Amaral and K. Irwin, arXiv 1802.04196(quant-ph)]. In this paper, previous work is encompassed by the use of torsion-free subgroups of Bianchi groups for deriving the quantum gate generators of uqc. A special role is played by a chain of Bianchi congruence $n$-cusped links starting with Thurston's link.

\end{abstract}

\maketitle

\vspace*{-.5cm}




\normalsize
\section{Introduction}

A Bianchi group $\Gamma_k=PSL(2,\mathcal{O}_k)<PSL(2,\mathbb{C})$ acts as a subset of orientation-preserving isometries of $3$-dimensional hyperbolic space $\mathbb{H}_3$, with $\mathcal{O}_k$ the ring of integers of the imaginary quadratic field $\mathcal{I}=\mathbb{Q}(\sqrt{-k})$. The quotient space $3$-orbifold $PSL(2,\mathcal{O}_k)\setminus \mathbb{H}_3$ has a set of cusps in bijection with the class group $\mathcal{I}$ \cite{MacLachlanBook, Thurston1997, AdamsBook}. A torsion-free subgroup $\Gamma_k(l)$ of index l of $\Gamma_k$ is the fundamental group $\pi_1$ of a $3$-manifold defined by a link such as the figure-eight knot [with $\Gamma_{-3}(12)$], the Whitehead link [with $\Gamma_{-1}(12)$] or the Borromean rings [with$\Gamma_{-1}(24)$]. The fundamental group of a knot or link complement (such as the complement  the figure-eight knot $K4a1$, the Whitehead link $L5a1$ or Borromean rings $L6a4$) was used to construct appropriate $d$-dimensional fiducial states for universal quantum compuring (uqc) \cite{PAAI2018}. The latter states come from the permutation structure of cosets in some $d$-coverings of subgroups $\Gamma_k(l)$ [alias the subgroups of index $d$ of subgroups $\Gamma_k(l)$ given a selected pair $(k,l)$].

In this paper, one starts by upgrading these models of uqc by using other torsion-free subgroups of Bianchi groups and the corresponding $3$-manifold such as the Berg\'e manifold that comes from the Berg\'e link $L6a2$ [with $\Gamma_{-3}(24)]$ or the so-called magic manifold that comes from the link $L6a5$  [with $\Gamma_{-7}(6)]$. The latter link is a small congruence link \cite{Baker2016} and belongs to a chain of eight links starting with Thurston's eight-cusped congruence link [with $\Gamma(-3)$ and ideal $\left\langle (5+\sqrt{-3})/2\right\rangle$] and ending with the Whitehead link and the figure-eight knot. In Sec. 2, it is presented the permutation based model of universal quantum computing developed by the authors and its relationship to some minimal informationally complete POVMs or MICs \cite{PlanatRukhsan,PlanatGedik,PlanatModular}. In Sec. 3, important small index torsion-free subgroups of Bianchi groups and the corresponding substructure of $3$-manifolds are derived. In Sec. 4, one specializes in the connection of the aforementioned Bianchi subgroups to our version of uqc. One focuses on a remarkable chain of $n$-cusped links, $n=8..1$ obtained thanks to $(\pm 1, 1)$-slope Dehn filling starting with congruence Thurston's link. Their possible role for uqc and the relevant $3$-manifolds is discussed.

\section{Minimal informationally complete POVMs with permutations}

\begin{figure}[h]
\centering 
\includegraphics[width=5cm]{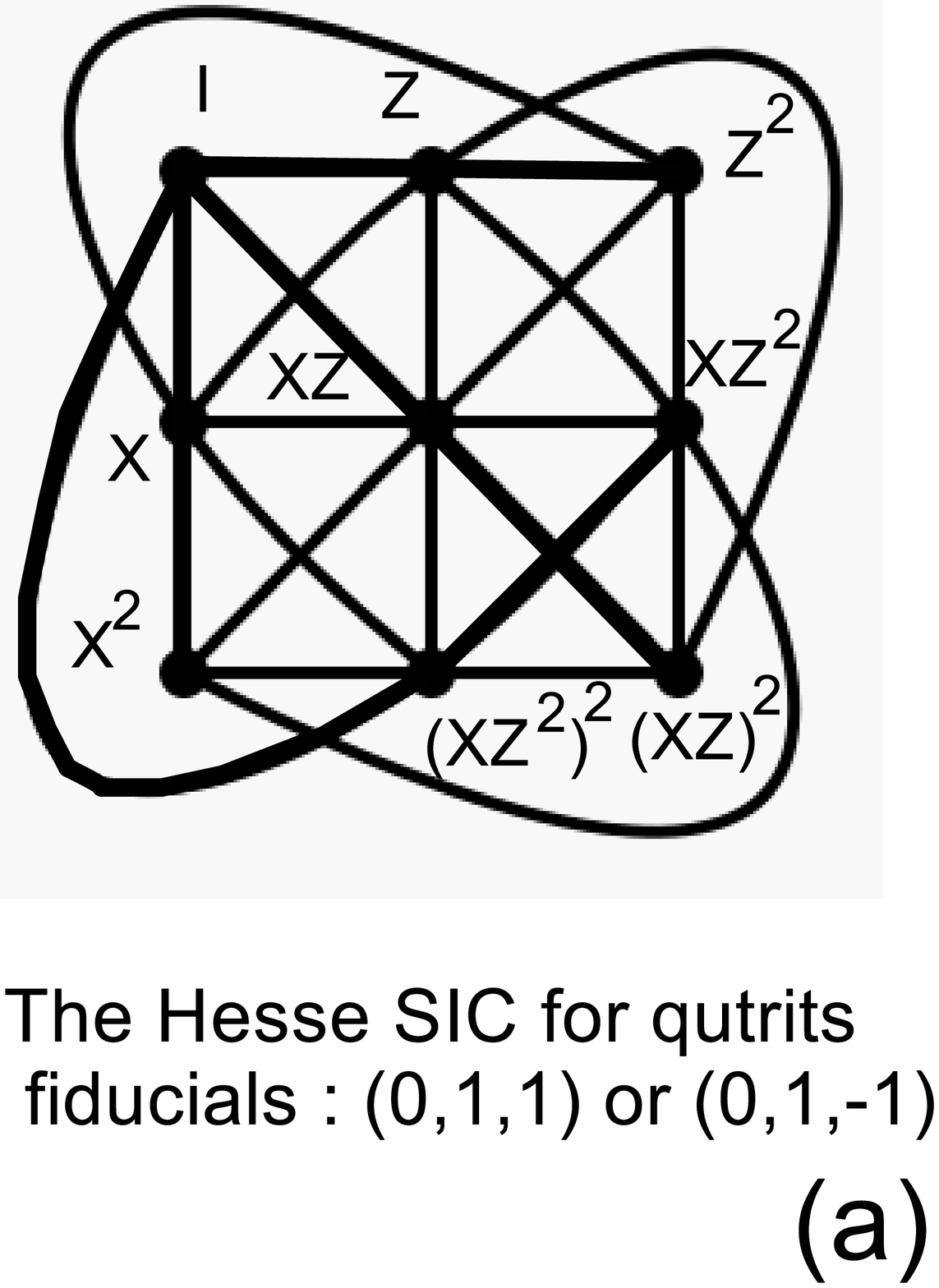}
\includegraphics[width=5cm]{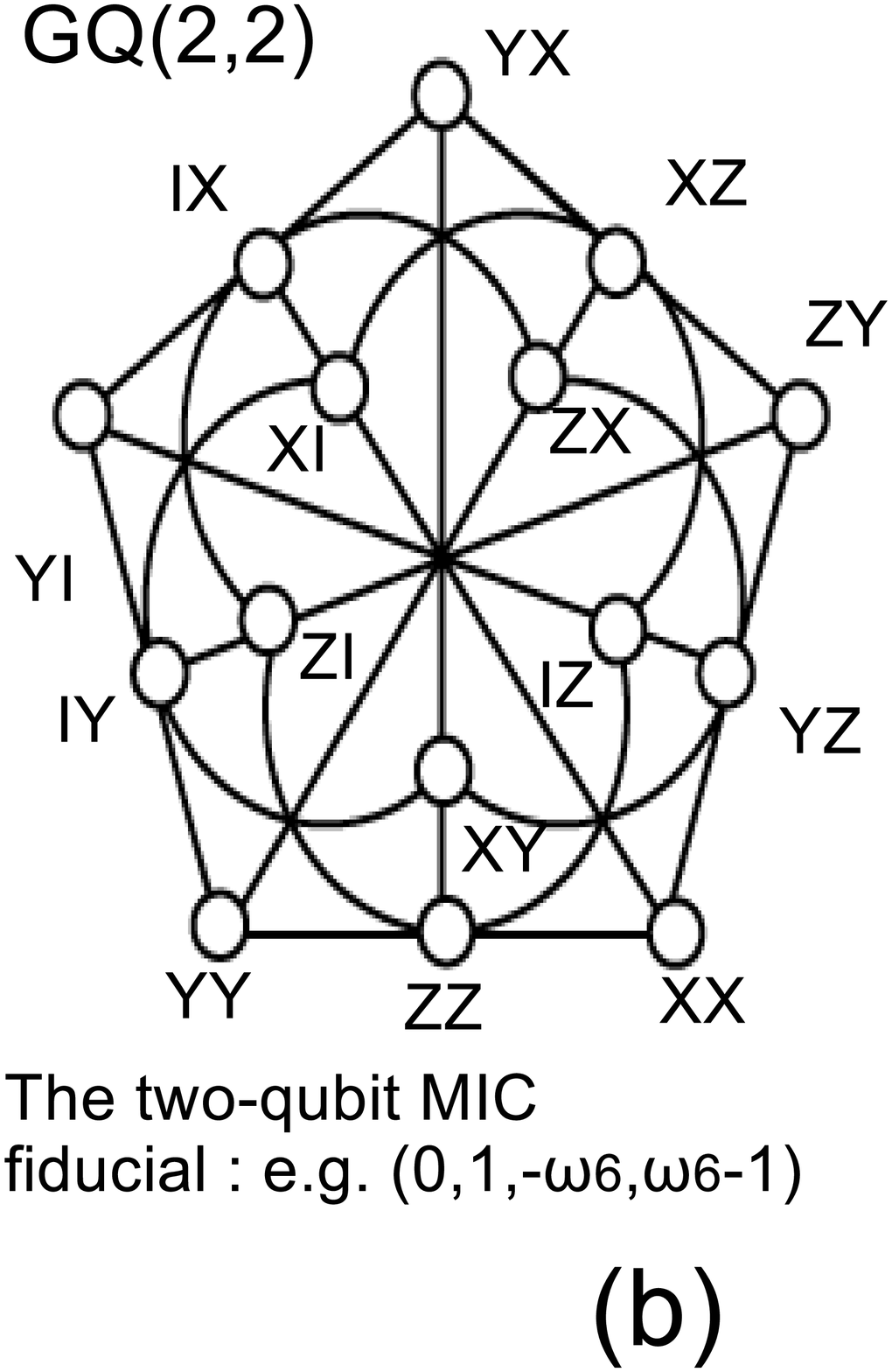}
\includegraphics[width=5cm]{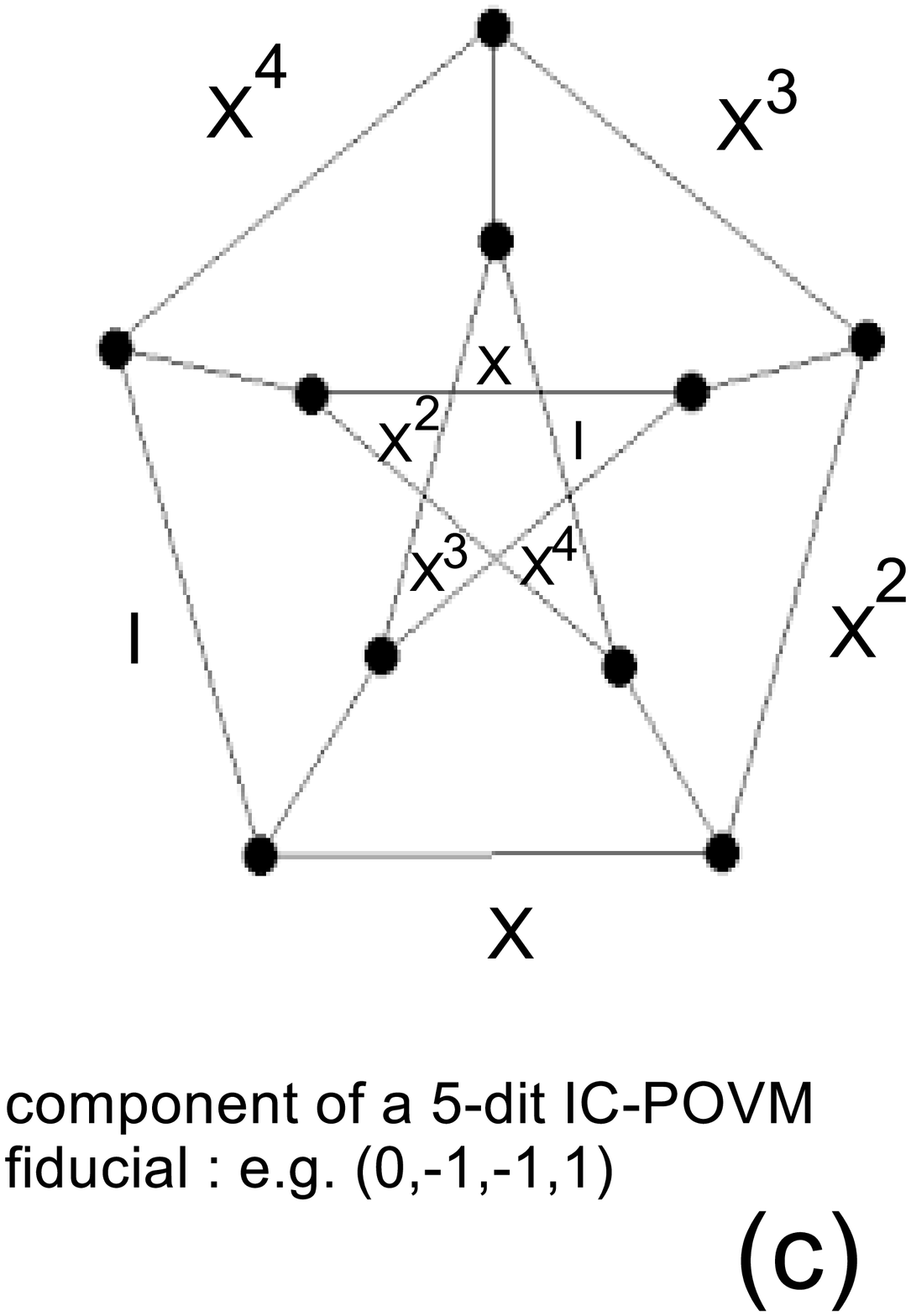}
\includegraphics[width=5cm]{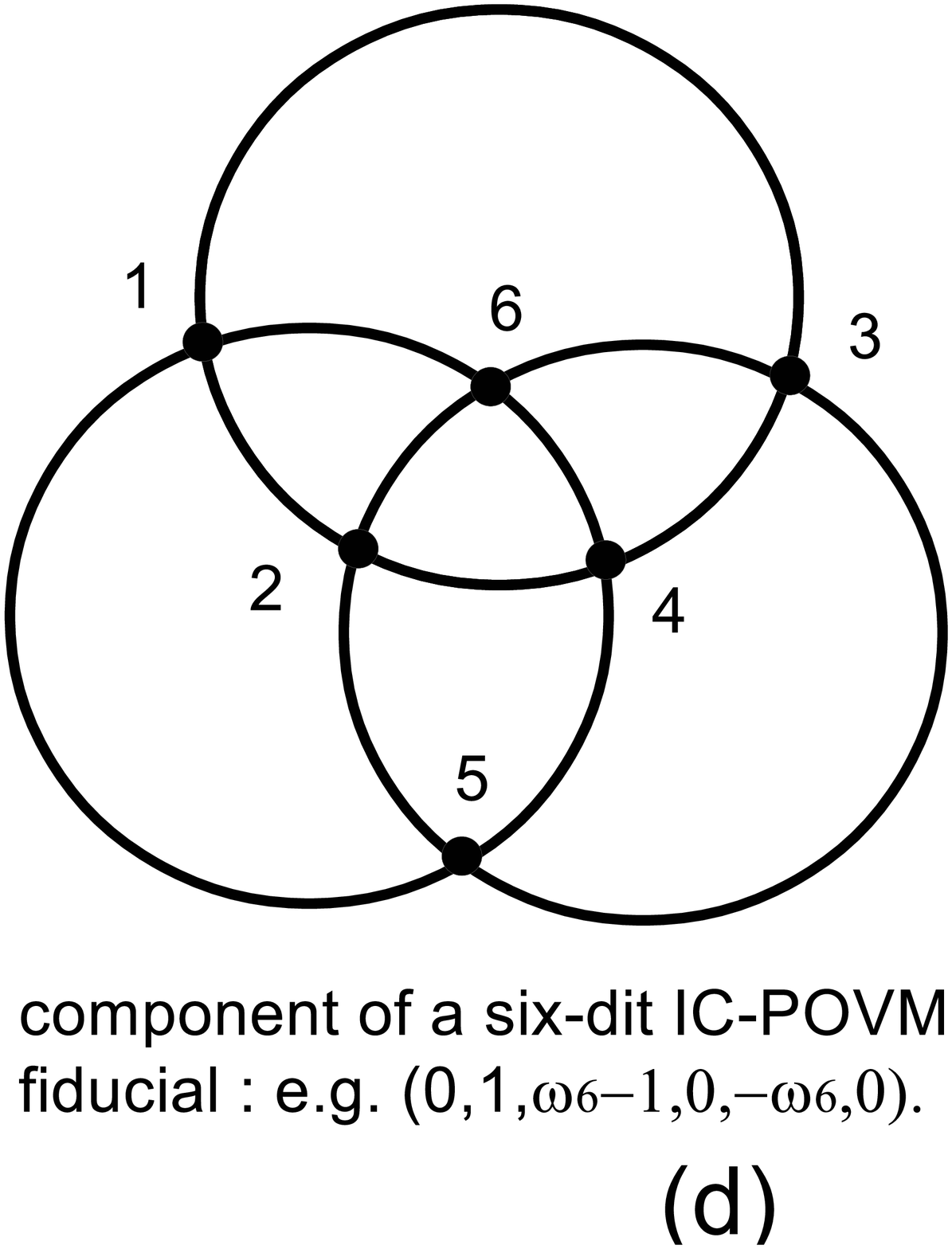}
\caption{Geometrical structure of low dimensional MICs: (a) the qutrit Hesse SIC, (b) the two-qubit MIC that is the generalized quadrangle of order two $GQ(2,2)$, (c) the basic component of the $5$-dit MIC that is the Petersen graph. The coordinates on each diagram are the $d$-dimensional Pauli operators that act on the fiducial state, as shown. For the $6$-dit case in (d), the coordinates are made explicit in the caption of \cite[Fig. 4]{PlanatModular}.}
\end{figure}

In our approach \cite{PlanatGedik, PlanatModular}, minimal informationally complete POVMs (MICs) are derived from appropriate fiducial states under the action of the (generalized) Pauli group. The fiducial states also allow to perform universal quantum computation \cite{PlanatRukhsan}.

A POVM is a collection of positive semi-definite operators $\{E_1,\ldots,E_m\}$ that sum to the identity. In the measurement of a state $\rho$, the $i$-th outcome is obtained with a probability given by the Born rule $p(i)=\mbox{tr}(\rho E_i)$. For a minimal IC-POVM, one needs $d^2$ one-dimensional projectors $\Pi_i=\left|\psi_i\right\rangle \left\langle \psi_i \right|$, with $\Pi_i=d E_i$, such that the rank of the Gram matrix with elements $\mbox{tr}(\Pi_i\Pi_j)$, is precisely $d^2$. A SIC-POVM (the $S$ means symmetric) obeys the relation
$\left |\left\langle \psi_i|\psi_j \right \rangle \right |^2=\mbox{tr}(\Pi_i\Pi_j)=\frac{d\delta_{ij}+1}{d+1},$
that allows the explicit recovery of the density matrix as in \cite[eq. (29)]{Fuchs2004}.

New MICs (i.e. whose rank of the Gram matrix is $d^2$) and with Hermitian angles $\left |\left\langle \psi_i|\psi_j \right \rangle \right |_{i \ne j} \in A=\{a_1,\ldots,a_l\}$ have been discovered \cite{PlanatModular}. A SIC is equiangular with $|A|=1$ and $a_1=\frac{1}{\sqrt{d+1}}$. The states encountered are considered to live in a cyclotomic field $\mathbb{F}=\mathbb{Q}[\exp(\frac{2i\pi}{n})]$, with $n=\mbox{GCD}(d,r)$, the greatest common divisor of $d$ and $r$, for some $r$. 
The Hermitian angle is defined as $\left |\left\langle \psi_i|\psi_j \right \rangle \right |_{i \ne j}=\left\|(\psi_i,\psi_j)\right\|^{\frac{1}{\mbox{\footnotesize deg}}}$, where $\left\|.\right\|$ means the field norm
 of the pair $(\psi_i,\psi_j)$ in $\mathbb{F}$ and $\mbox{deg}$ is the degree of the extension $\mathbb{F}$ over the rational field $\mathbb{Q}$ \cite{PlanatGedik}.

The fiducial states for SIC-POVMs are quite difficult to derive and seem to follow from algebraic number theory \cite{Appleby2017}. Except for $d=3$, the MICs derived from permutation groups are not symmetric and most of them can be recovered thanks to subgroups of index $d$ of the modular group $\Gamma$ \cite[Table 1]{PlanatModular}. For instance, for $d=3$, the action of a magic state of type $(0,1,\pm 1)$ results into the Hesse SIC shown in Fig 1a, arising from the congruence subgroup $\Gamma_0(2)$ of $\Gamma$. For $d=4$, the action of the two-qubit Pauli group on the magic/fiducial state of type $(0,1,-\omega_6,\omega_6-1)$ with $\omega_6=\exp(\frac{2i\pi}{6})$ results into a MIC whose geometry of triple products of projectors $\Pi_i$, arising from the congruence subgroup $\Gamma_0(3)$ of $\Gamma$, turns out to correspond to the commutation graph of Pauli operators, see Fig. 1b and \cite[Fig. 2c]{PlanatModular}. For $d=5$, the congruence subgroup  $5A^0$ of $\Gamma$ is used to get a MIC whose geometry consists of copies of the Petersen graph (see Fig. 1c and \cite[Fig. 3c]{PlanatModular}. For $d=6$, all five congruence subgroups $\Gamma'$, $\Gamma(2)$, $3C^0$, $\Gamma_0(4)$ or $\Gamma_0(5)$ point out the geometry of Borromean rings (see \cite[Fig. 4c]{PlanatModular} and Table 1 of \cite{PAAI2018}).

The modular group $\Gamma$ first served as a motivation for investigating the trefoil knot manifold $3_1$ in relation to uqc and the corresponding MICs, then the uqc problem was put in the wider frame of Poincar\'e conjecture, the Thurston's geometrization conjecture and the related $3$-manifolds.  E.g. MICs may follow from hyperbolic or Seifert $3$-manifolds as shown in Tables 2 to 5 of \cite{PAAI2018}. A further step is obtained here by restricting our choice of $3$-manifolds to low degree coverings of Bianchi subgroups. This is a quite natural procedure to base our uqc models on Bianchi groups $\Gamma_k$ defined over the ring $\mathcal{O}_k$ of imaginary quadratic integers \cite{Grunewald1993} since they are a natural generalization of the modular group $\Gamma$.

\section{Torsion-free subgroups of Bianchi groups}

\begin{table}[h]
\begin{center}
\begin{tabular}{|c|l|c|c||r|r|}
\hline 
\hline
k & index& cusps & names& other name & note\\
\hline
 -1        &$12$          &$2$        &{\bf L5a1}, ooct$01_{00001}$        & \small{Whitehead link: WL}&L6a5(1,1)\\
.        &.          &.        &L13n5885, ooct$01_{00000}$        & sister of WL  &\scriptsize {L10n113(1,1)(1,1)(1,1)}\\
.        &$24$          &$3$        &{\bf L6a4}, ooct$02_{00005}$        &Borromean rings&\\
.        &.          &.        &L8n7, ooct$02_{00001}$        & &L5a1(1,1)\\
.        &.          &.        &L10n84, ooct$02_{00002}$        &  &\\
\hline
 -2        &$12$          &$2$        &L9a32, $9_{40}^2$        & &\\
.       &.          &.        &L9a33, $9_{24}^2$        & &\\
\hline
 -3        &$12$          &$2$        &{\bf K4a1}, otet$02_{00001}$        & figure-eight knot  &L5a1(1,1)\\
.        &.          &$1$        &m003, otet$02_{00000}$        & sister of K4a1 &\\
.       &$24$          &$2$        &{\bf L6a2}, $6_2^2$, otet$04_{00001}$        & Berg\'e manifold&\\
.       &.         &$1$        &m206, otet$04_{00002}$        & &\\
.       &.          &.        &m207, otet$04_{00001}$        & &\\
\hline
-7        &$6$          &$2$        &{\bf L6a1}, $6_3^2$        &  &\\
.        &.          &$3$        &{\bf L6a5}, $6_1^3$        & magic manifold  & congruence\\
.        &$12$         &$3$        &L10n81, $10_{15}^3$        &  &\scriptsize { cyc. $2$-cover of L6a5}\\
.        &.          &$4$        &  L12n2205      &   & \\
\hline
-15        &$6$         &$6$      &  $\left\langle (1+\sqrt{-15})/2\right\rangle$    & Baker\&Reid \cite{Baker2016} &congruence\\
\hline
\hline
\end{tabular}
\caption{Low index torsion-free subgroups of Bianchi groups $\Gamma_k=PSL(2,\mathcal{O}_k)$, $k=-1,-2,-3,-7$, their names in SnapPy in column 4 and their popular name in column 5, when known. The last column shows a connection to Dehn filling of slope $(1,1)$ of the corresponding manifold. Bold characters point out manifolds investigated in this paper  for their connection to MICs.} 
\end{center}
\end{table}

Table 1 provides a short list of low index torsion-free subgroups of Bianchi groups $\Gamma_k$, $k \in \{-1,-2,-3,-7,-15\}$ \cite{Grunewald1993}, obtained from the software Magma \cite{Magma} and their identification as fundamental groups of $3$-manifolds obtained from the software SnapPy \cite{SnapPy}. In particular, one recovers the Whitehead link and the Borromean rings ($k=-1$ with index $l=12$ and $24$, respectively) as well as the figure-eight knot ($k=-3$ with index $12$) whose relationship to uqc was explored in \cite{PAAI2018}. Other links of importance in our paper are the link L6a1, the Berg\'e link L6a2 and the magic link L6a5 whose connections to uqc are summarized in Table 3.

According to \cite{Mednykh2006}, a $d$-fold covering of the fundamental group $\pi_1(M)$ of a $3$-manifold $M$ is uniquely determined by the conjugacy class of a subgroup of index $d$ of $\pi_1$. To recognize the $d$-covering from the conjugacy class and conversely, one makes use of the cardinality signature of subgroups of $\pi_1(M)$, denoted $\eta_d(\pi_1(M))$, $d=1..d_{\mbox{max}}$ and one identifies $\eta_d(\pi_1(M))$ both in a particular $d$-covering of $M$ (with SnapPy) and in the representative of a particular conjugacy class of subgroups of $\pi_1(M)$ (with Magma).

Table 2 provides the identification of all $3$-manifolds coming from $d$-fold coverings ($d \le 5$) of the Berg\'e link and identifies the ones related to MICs. Observe that it becomes a cumbersome task for manifolds with many cusps. E.g., for the manifold otet$16_{00025}$ that one finds corresponding to the two-qubit MIC, one gets $\eta_d(M)=[15,70,642,2206,30192,\cdots]$ when $d =[2,3,4,5,6,\cdots]$.

\begin{table}[h]
\begin{center}
\begin{tabular}{|c|l|c|c|r|r|}
\hline 
\hline
d & homology&  \footnotesize{cusps} & sym & names & note\\
 \hline
2 &cyc: $1+1$   &$2$ &$D_4$&  otet$08_{00002}$, $L10n46$  &\\
. &cyc: $\frac{1}{5}+1+1$ &. &$G_{16}$ &  otet$08_{00010}$   &\\
\hline
3 &cyc: $1+1+1+1$   &$4$ &$\frac{Z}{2}\times O$&  otet$12_{00009}$,  $L12n2208$ &\\
. &cyc: $\frac{1}{3}+\frac{1}{3}+1+1$   &. &$D_6$ &  otet$12_{00020}$ &\\
\hline
4 &cyc: $1+1+1$   &$2$ & $D_8$&  otet$16_{00013}$, $L14n17878$  &\\
. &cyc: $\frac{1}{3}^{+2}+1^{+2}$   &. & $G_{16}$&  otet$16_{00091}$   &\\
. &cyc: $\frac{1}{5}+1+1$   &. & $\frac{Z}{2}\times D_4$&  otet$16_{00063}$   &\\
. &irr: $1+1+1+1$   &$4$ & $\frac{Z}{2}+\frac{Z}{2}$&  otet$16_{00025}$  & 2QB MIC\\
. &reg: $\frac{1}{5}+1+1$   &$2$ & $G_{32}$&  otet$16_{00092}$  & \\
\hline
5 &cyc: $1+1$   &2 & $D_{10}$&  otet$20_{00443}$,  &\\
. &cyc: $\frac{1}{2}^{+4}+1+1$   &. & $G_{20}$&  otet$20_{01343}$,  &\\
. &irr: $1^{+4}$   &$4$ &$\frac{Z}{2}$ & \small{ otet$20_{00543}$; otet$20_{00041}$; otet$20_{00549}$}  &$5$-dit MIC\\
. &irr: $\frac{1}{2}+1^{+4}$   &$.$ &. &  otet$20_{00574}$  &$5$-dit MIC\\
. &irr: $\frac{1}{3}+\frac{1}{3}+1+1$   &$2$ &$D_4$ &  otet$20_{00665}$  &$5$-dit MIC\\
. &irr: $1^{+6}$   & $6$ &$\frac{Z}{2}\times D_4$ &  otet$20_{00573}$  &\\
\hline
\hline
\end{tabular}
\caption{Table of $3$-manifolds (in column 5) found from subgroups of finite index $d$ of the fundamental group of Berg\'e link (alias the $d$-fold coverings of $L6a2$). The terminology is that of SnapPy \cite{SnapPy} with a first homology group such as $\mathbb{Z}$ denoted as $1$ in column 2 or $1^{+s}$ for $\mathbb{Z}+\cdots \mathbb{Z}$ ($s$ times). Column 4 is the symmetry group of the manifold with $O$ the octahedral group and $G_a$ a group of order $a$. The three manifolds separated by a semicolon are distinct although with the same homology, symmetry and number of cusps.} 
\end{center}
\end{table}

\section{A Bianchi factory for quantum computing}

First one provides a reminder about the concept of Dehn filling that will be useful in Sec. 4.2. Let us start with a Lens space $L(p,q)$ that is $3$-manifold obtained by gluing the boundaries of two solid tori together, so that the meridian of the first solid torus goes to a $(p,q)$-curve on the second solid torus [where a $(p,q)$-curve wraps around the longitude $p$ times and around the meridian $q$ times]. Then we generalize this concept to a knot exterior, i.e. the complement of an open solid torus knotted like the knot. One glues a solid torus so that its meridian curve goes to a $(p,q)$-curve on the torus boundary of the knot exterior, an operation called Dehn surgery \cite[p. 275]{Thurston1997},\cite[p 259]{AdamsBook}, \cite{Gordon1998},\cite{Martelli2014}. According to Lickorish's theorem, every closed, orientable, connected 3-manifold is obtained by performing Dehn surgery on a link in the 3-sphere.

Another operation useful for creating new $3$-manifolds from old ones is that of drilling out along a closed geodesic. Let us quote Thurston \cite{ThurstonLink}
{\it In SnapPea there is a module for drilling out geodesic loops that are identified geometrically. Given a hyperbolic manifold, the program will present a list of loops up to a specified combinatorial complexity sorted in order of hyperbolic length. SnapPea can remove any of these geodesics from the manifold, in essence deforming its cone angle to $0$ so as to obtain a complete hyperbolic structure for its complement... Repeatedly drilling out a relatively short geodesics always seems to start crystallizing a $3$-manifold, so that its cusp neighbourhoods can be adjusted in size to pack together in regular crystalline patterns, most commonly in a square- looking pattern as for the Borromean rings}.

One shows below how the Bianchi subgroups are related to uqc and the related MICs and how useful is Dehn filling and conversely drilling in our context.

\begin{table}[h]
\begin{center}
\begin{tabular}{|c|l|c|c|r|r|}
\hline 
\hline
source & index& cusps & homology& names or $\eta_d(M)$ & uqc\\
\hline
L6a1 & $3$ & $3$ &irr: $1^{+4}$ & L12n2181 & Hesse SIC\\
. & $4$ & $5$ &irr: $1^{+5}$ & L14n63905 & 2QB MIC\\
\hline
Berg\'e L6a2 & $4$ & $4$ &irr: $1^{+4}$                 & otet$16_{00025}$ & 2QB MIC\\
.            & $5$ & $4$ &irr: $1^{+4}$                 & \scriptsize{otet$20_{00543}$; otet$20_{00041}$; otet$20_{00549}$}  & $5$-dit MIC\\
.            & $5$ & $2$ &irr: $\frac{1}{2}+1^{+4}$     & otet$20_{00574}$ & .\\
.            & $5$ & $2$ &irr: $\frac{1}{3}^{+2}+1^{2}$ & otet$20_{00665}$ & .\\
\hline
magic L6a5 & $3$ & $5$ &irr: $1^{+5}$ & L14n63788 & Hesse SIC\\
.          & $4$ & $4$ &irr: $\frac{1}{2}+1^{+4}$ &$[31,174,4324,82357,\cdots]$ & 2QB MIC\\
\hline
Borr L6a4 & $3$ & $4$ &irr: $1^{+4}$ & ooct$06_{00466}$ & Hesse SIC\\
.         & $3$ & $4$ &irr: $1^{+5}$ &  ooct$06_{00398}$&.\\
.         & $4$ & $4$ &irr: $\frac{1}{2}^{+2}+1^{+4}$ &$[63,300,10747,\cdots]$  & 2QB MIC\\
.         & $4$ & $6$ &irr: $\frac{1}{2}+1^{+6}$ &$[127,2871,478956,\cdots]$ & .\\
\hline
$4$-link L8n7 & $3$ & $6$ &irr: $1^{+6}$ & ooct$06_{00340}$ & Hesse SIC\\
\hline
$5$-link L10n113 & $3$ & $8$ &irr: $1^{+8}$ & ocube$06_{03005}$;ocube$06_{03322}$ & Hesse SIC\\
\hline
$6$-link L12n2256 & $3$ & $9$ &irr: $1^{+9}$ &$[511,20122,\cdots]$ & Hesse SIC\\
\hline
\hline
\end{tabular}
\caption{A few $3$-manifolds (in column 5) occuring as a model of uqc (in column 6). When not identified in SnapPy the cardinality structure $\eta_d(M)$ of finite index subgroups of the $3$-manifold $M$ is written explicitely. The fourth column identifies the first homology type, with $1$ meaning the ring of integers $\mathbb{Z}$. The $3$-manifolds occuring as a model of uqc with source the trefoil knot K3a$=3_1$, the figure-eight knot K4a1 and the Whitehead link L5a1 have been investigated in \cite{PAAI2018}.} 
\end{center}
\end{table}

\begin{figure}[h]
\centering 
\includegraphics[width=6cm]{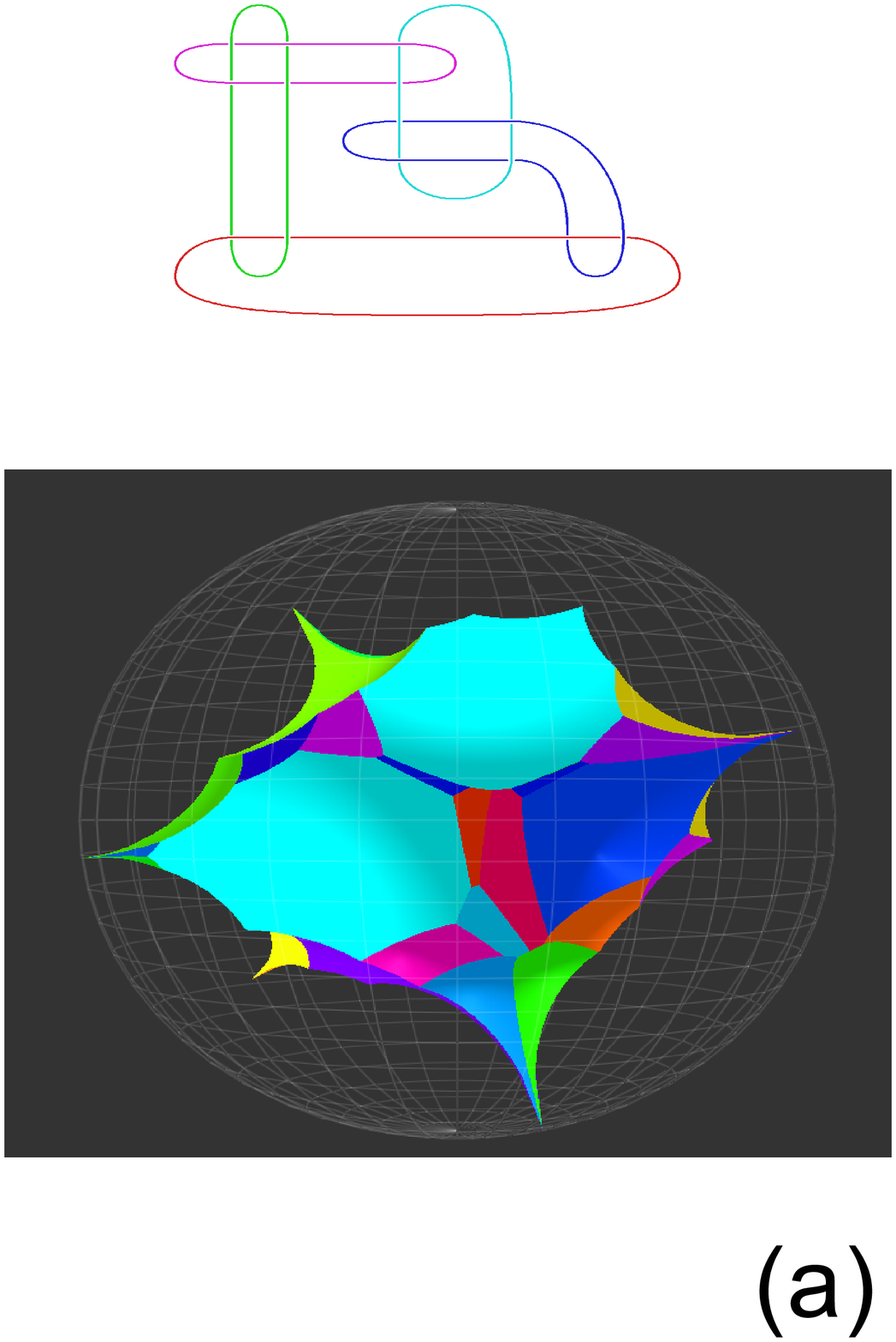}
\includegraphics[width=6cm]{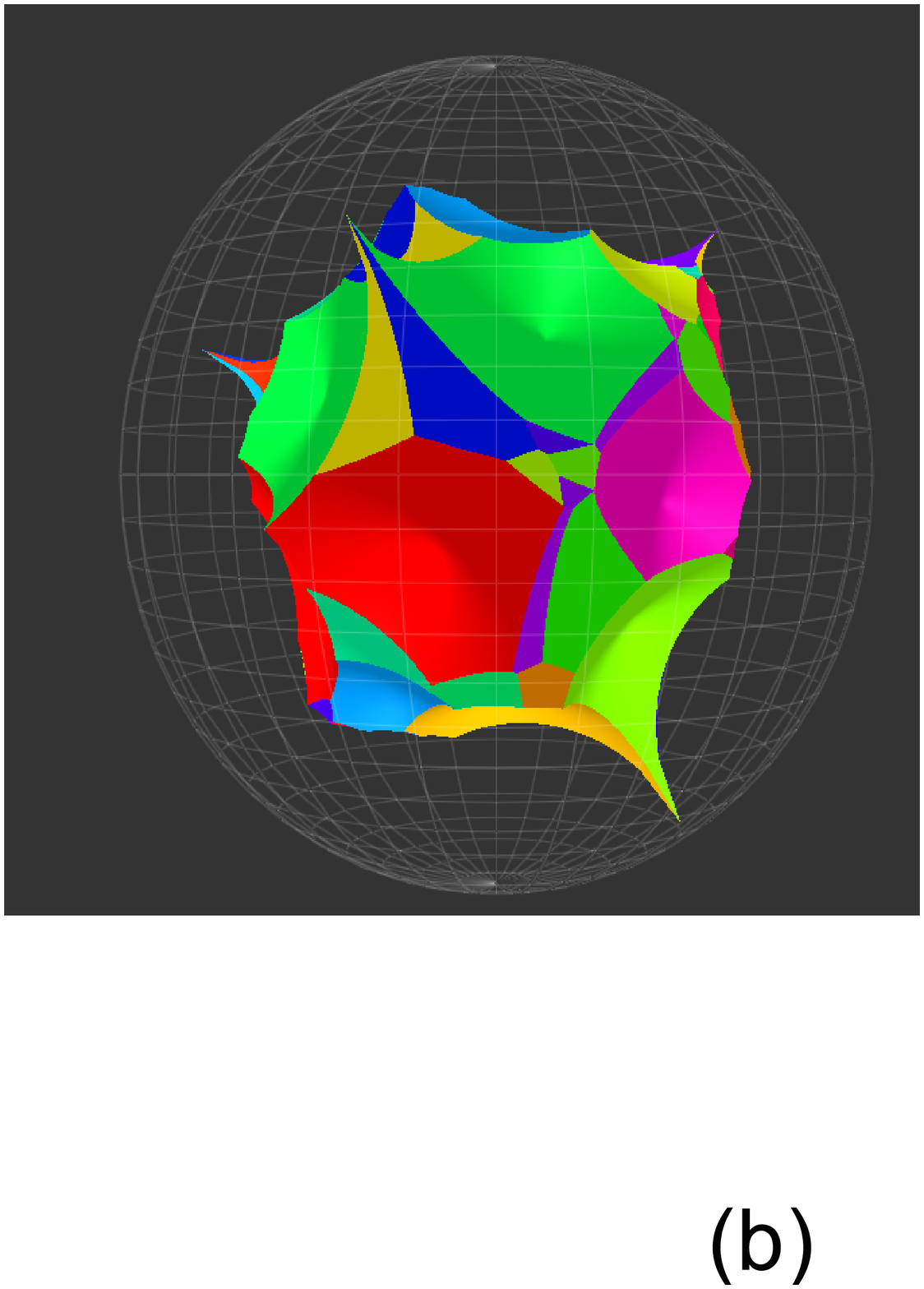}
\caption{(a) The $3$-manifold  $L14n63788$ corresponding to the $3$-fold covering of $L6a5$ relevant to qutrit universal quantum computing, (b) The $3$-manifold  corresponding to the $4$-fold covering of $L6a5$ relevant to two-qubit universal quantum computing.
 }
\end{figure} 

\subsection{Universal quantum computing from the Bianchi factory}

In this subsection, one specializes on uqc based on the qutrit Hesse SIC (shown in Fig. 1a) and on the two-qubit geometry of the generalized quadrangle $GQ(2,2)$ ( shown in Fig. 1b) obtained with the \lq magic' link L6a5. The qutrit uqc follows from a link called L14n63788 (with Poincar\'e polyhedron shown in Fig. 2a) and the two-qubit uqc follows from a $3$-manifold with Poincar\'e polyhedron (shown in Fig. 2b). The former case corresponds to the $3$-fold irregular covering of the fundamental group $\pi_1(L6a5)$ whose $3$-manifold has first homology $\mathbb{Z}^{\oplus 5}$, $5$ cusps and volume $\approx 16.00$. The latter case corresponds to the $4$-fold irregular covering whose $3$-manifold has first homology $\frac{\mathbb{Z}}{2}\oplus \mathbb{Z}^{\oplus 4}$, $4$ cusps and volume $\approx 21.34$. 

For the other links investigated in Table 3 one proceeds in the same way.

\subsection{Congruence links in the Bianchi factory}

\begin{figure}[ht]
%
\centering 
\includegraphics[width=3cm]{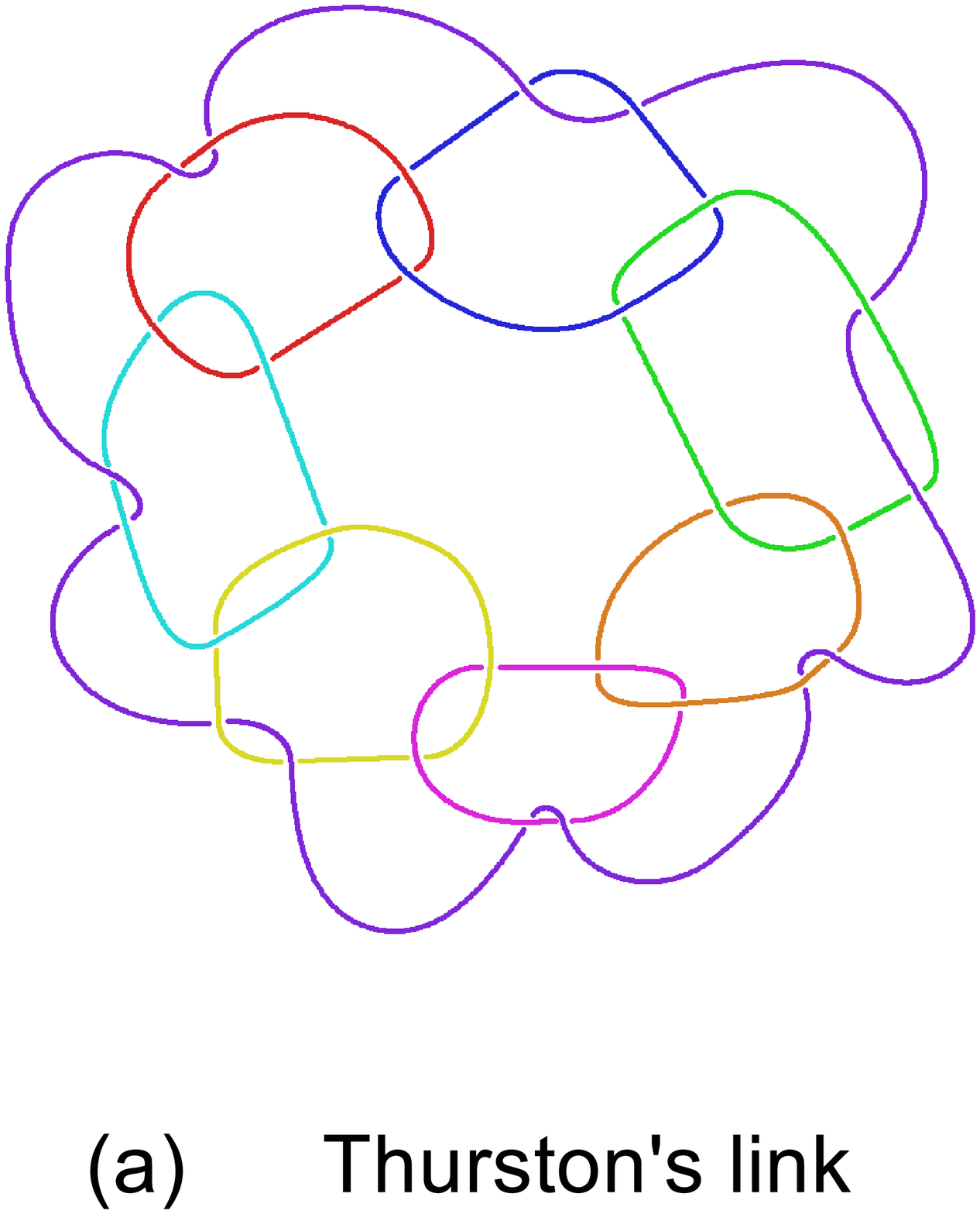}
\includegraphics[width=3cm]{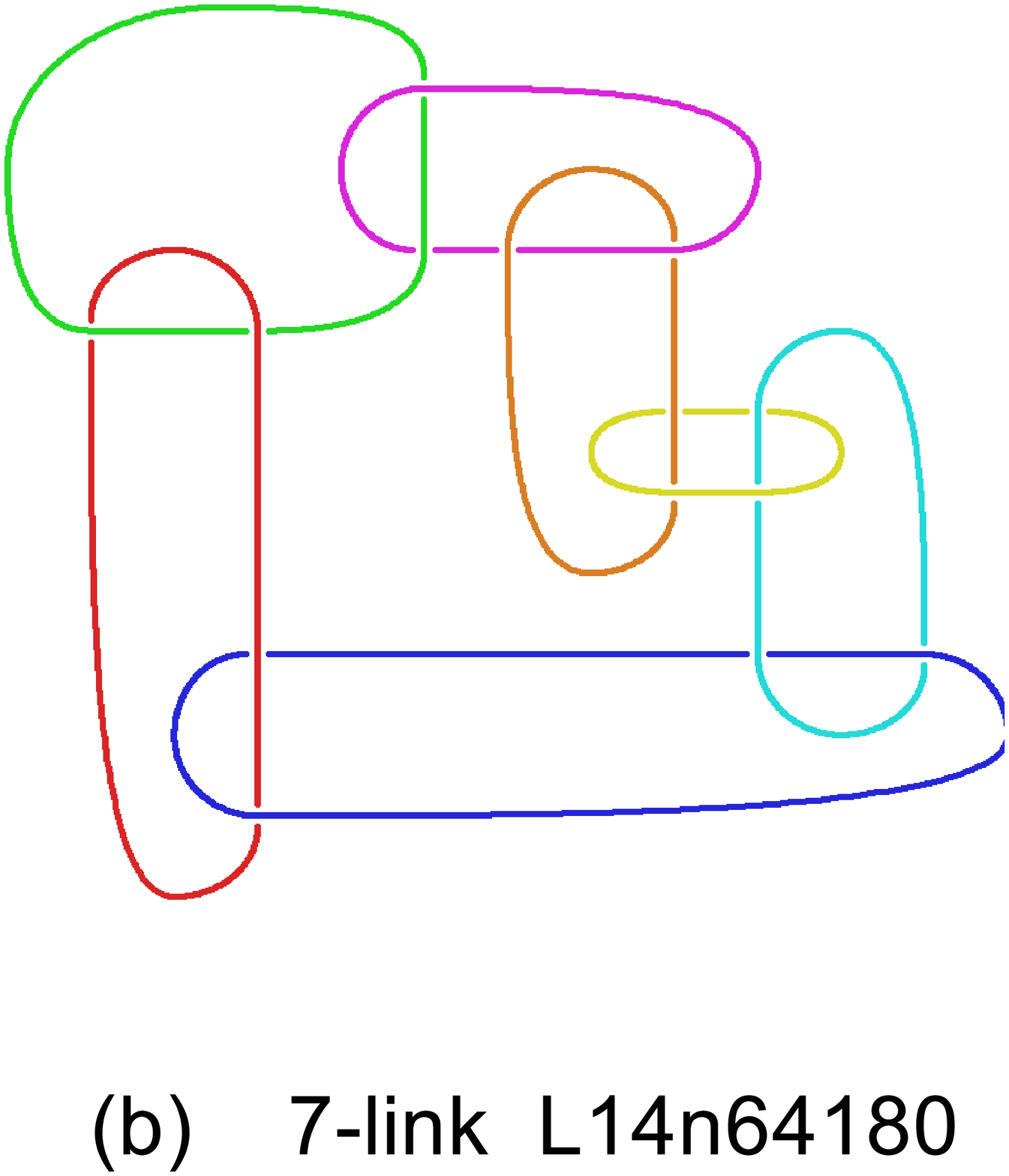}
\includegraphics[width=3cm]{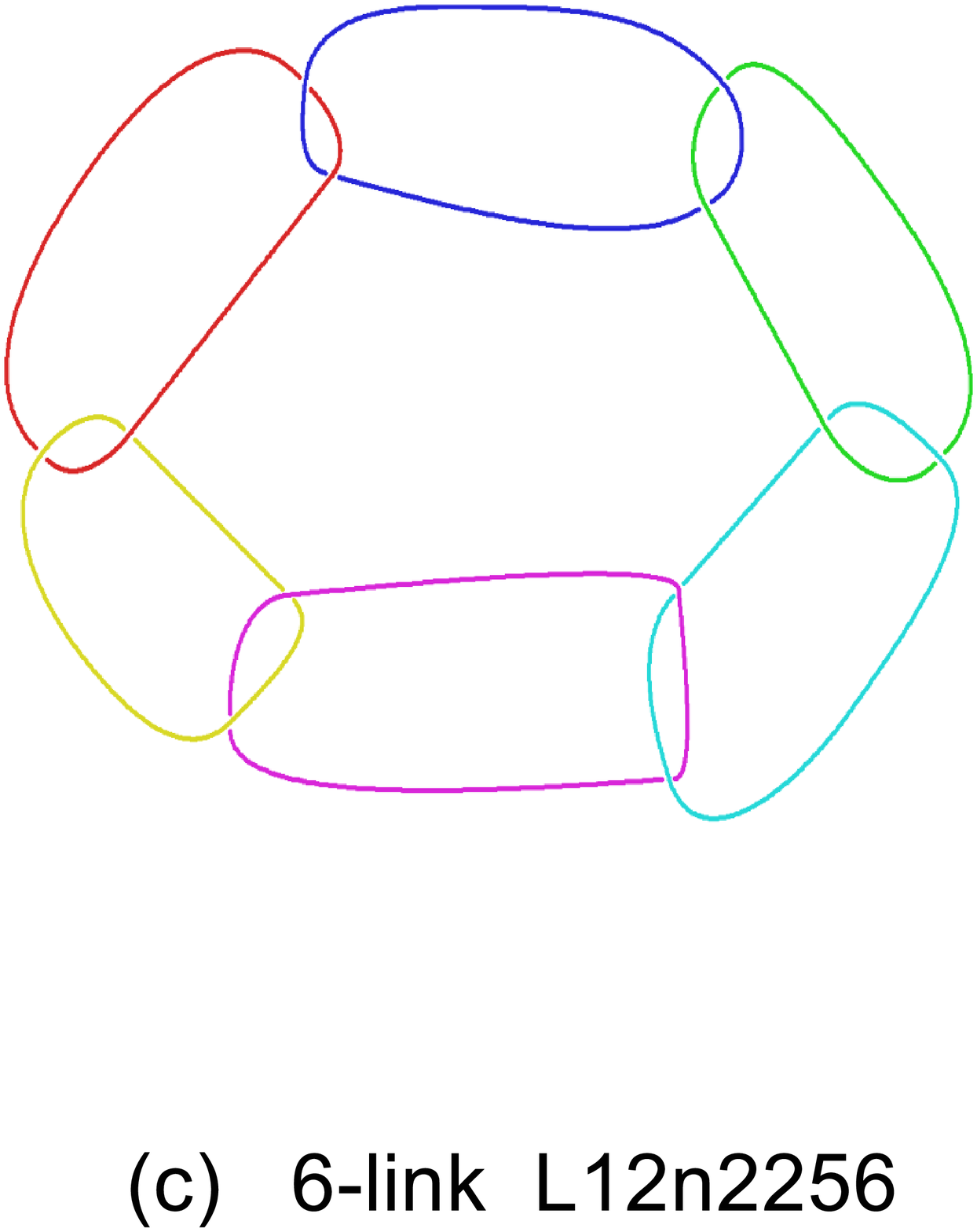}
\includegraphics[width=3cm]{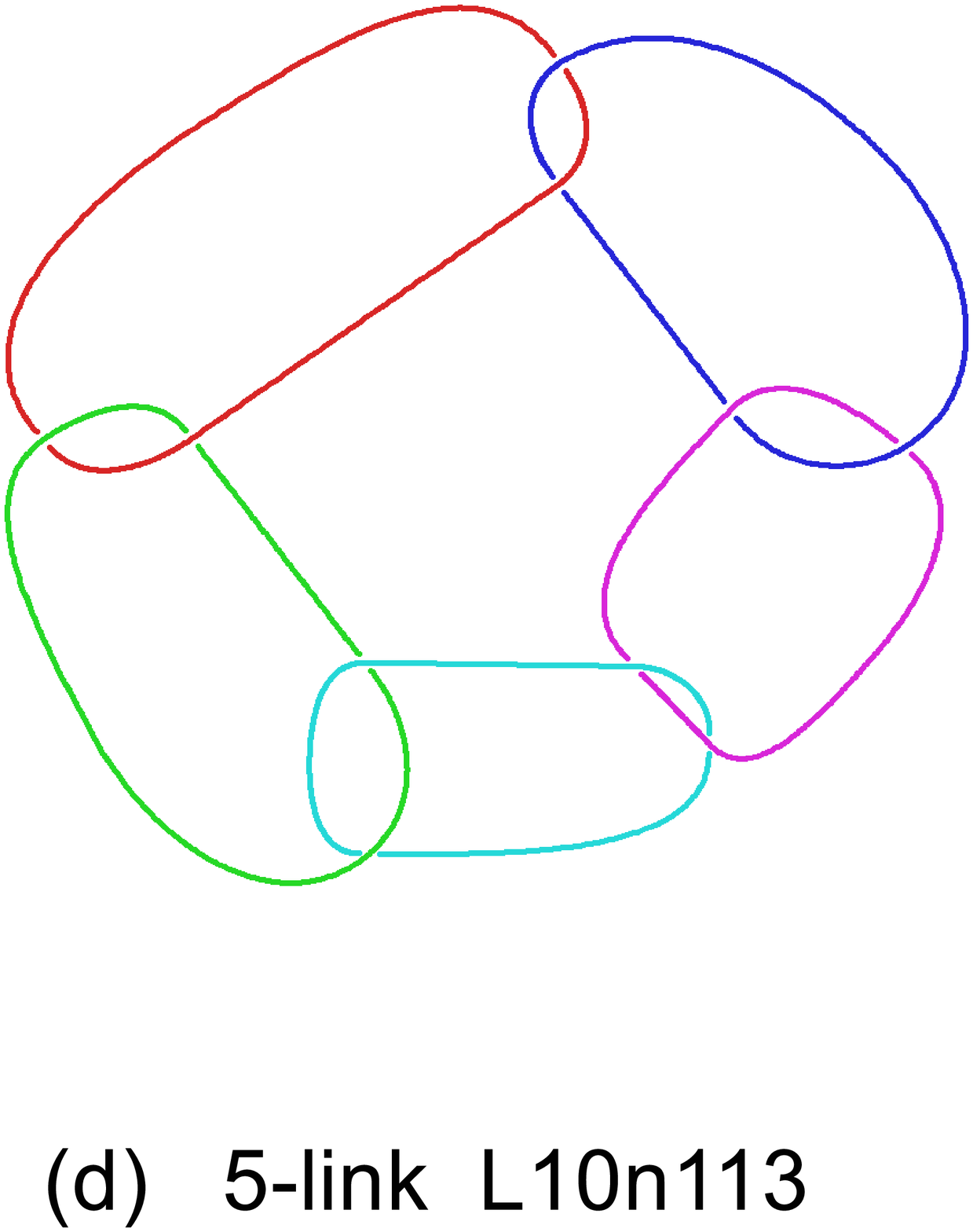}
\includegraphics[width=3cm]{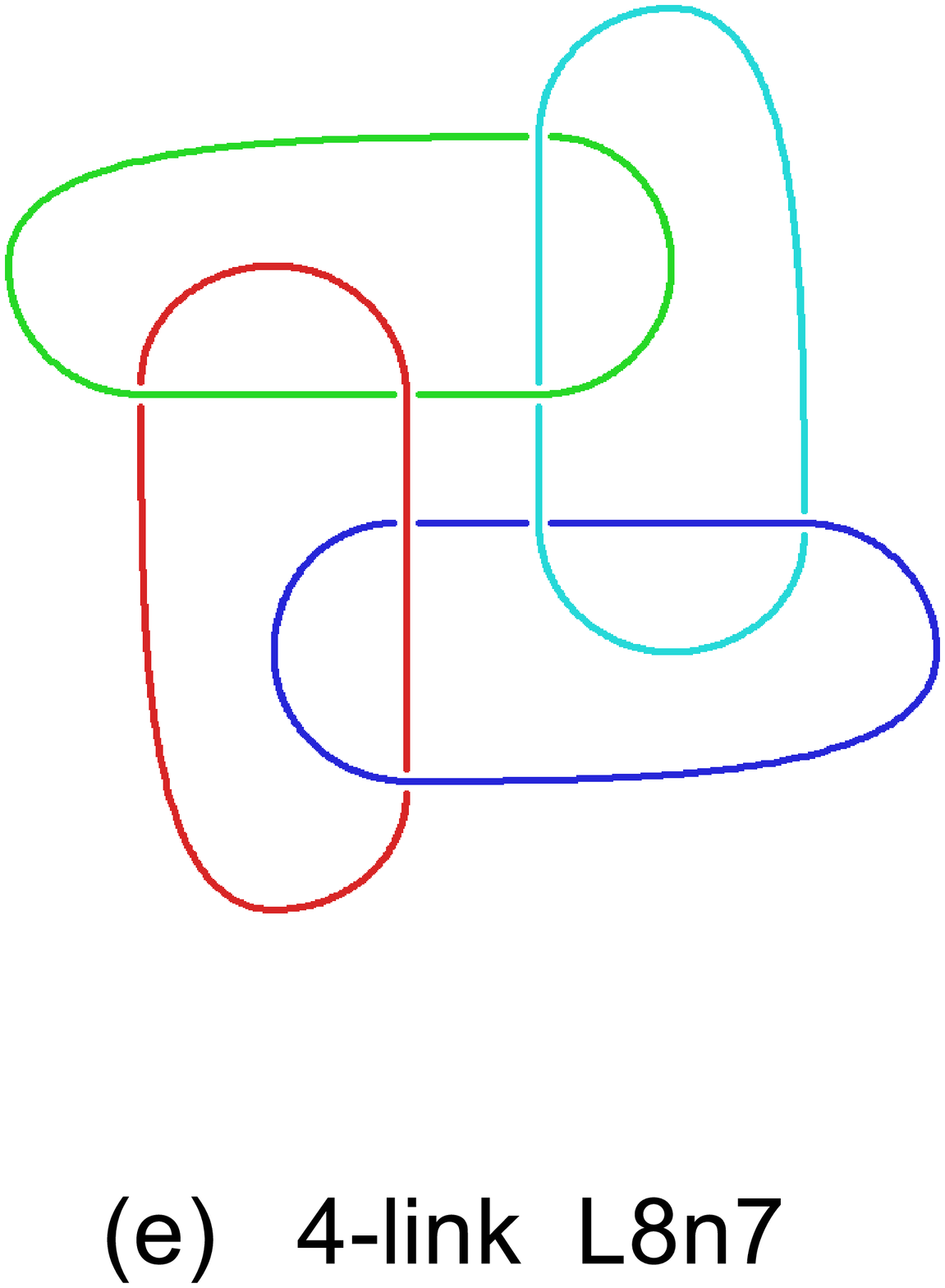}
\includegraphics[width=3cm]{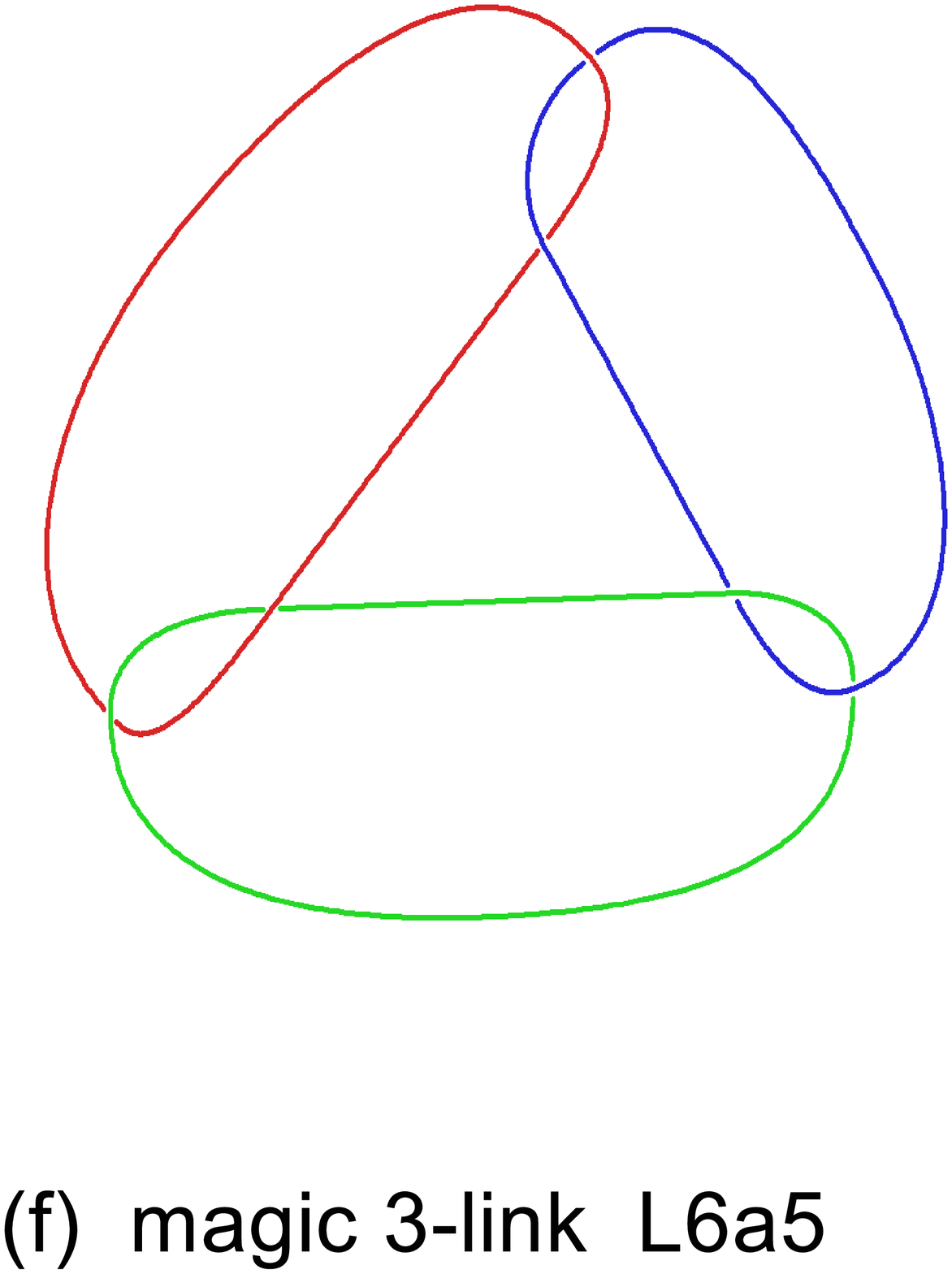}
\includegraphics[width=3cm]{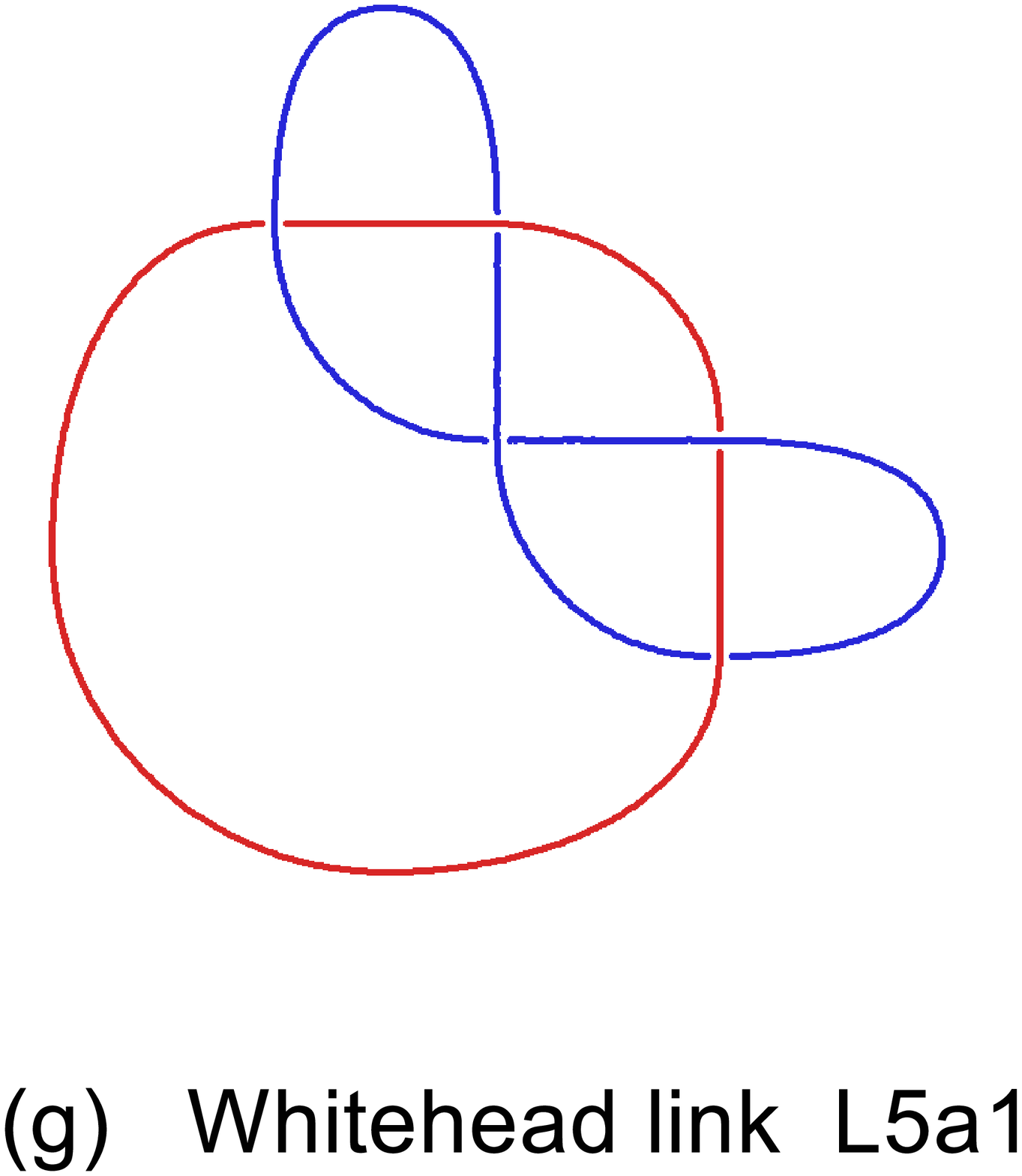}
\includegraphics[width=3cm]{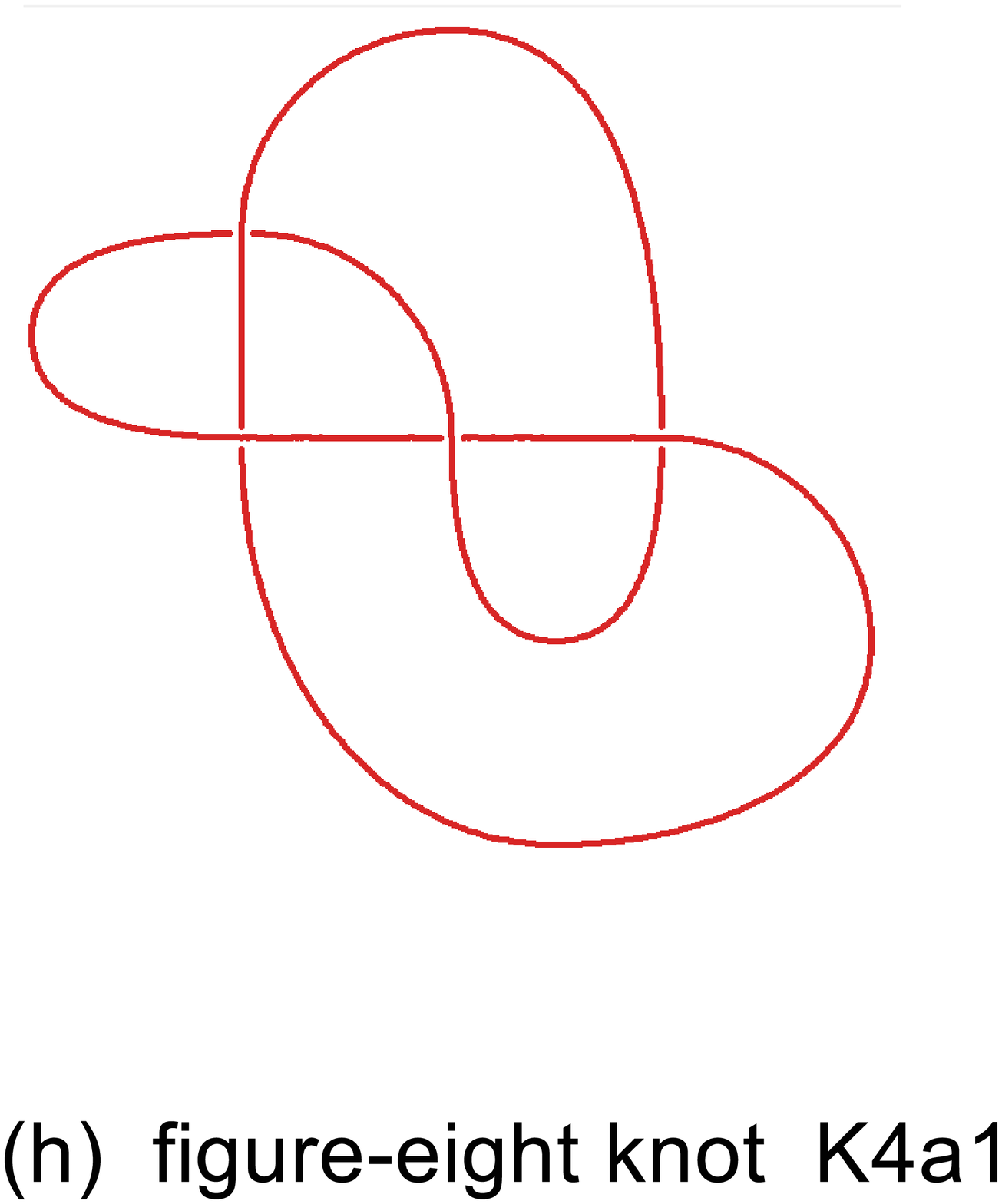}
\caption{The consecutive links from (a) to (h) obtained from (1,1) Dehn filling at a single cusp [except for $(-1,1)$ filling at  (e)]. The links (a), (c), (d), (f) are congruence. The sequence (d) to (h) already appears in \cite{Martelli2014}. Conversely, drilling out along the shortest closed geodesic (provided by SnapPy), one gets the chain $(h)\rightarrow (g)\rightarrow \cdots \rightarrow (b)$.}
\end{figure} 

As announced earlier, there is an interesting sequence of links starting at the Thurston's link \cite{ThurstonLink} and ending at the figure-eight knot that one obtains by applying $(\pm 1,1)$-slope Dehn fillings. The Dehn fillings of the last five links called $M_i$, $i=5..1$ were studied in \cite{Martelli2014}. Observe that $M_3$ is the magic link L6a5, $M_4$ is the $4$-link L8n7 and $M_5$ is the $5$-link L10n113 of Table 3. The full sequence is as follows

\begin{eqnarray}
 \mbox{{\bf Thurston's}}~\mbox{{\bf link}}:~~\left\langle(5+\sqrt{-3})/2  \right\rangle \stackrel{(1,1)}{\longrightarrow} 7-\mbox{link}~ L14n64180 \nonumber \\
 \stackrel{(1,1)}{\longrightarrow} \left\langle 2+0\sqrt{-1}  \right\rangle:~6~\mbox{link}~ L12n2256  \nonumber \\
 \stackrel{(1,1)}{\longrightarrow}\left\langle 2+0\sqrt{-3}  \right\rangle:~ 5-\mbox{link}~ L10n113 \nonumber \\
 \stackrel{(1,1)}{\longrightarrow} 4-\mbox{link}~ L8n7 \nonumber \\ \stackrel{(-1,1)}{\longrightarrow} ~\left\langle(1+\sqrt{-7})/2  \right\rangle:~ \mbox{magic}~ 3-\mbox{link}~L6a5\nonumber \\
\stackrel{(1,1)}{\longrightarrow} \mbox{Whitehead}~\mbox{link}~L5a1\stackrel{(1,1)}{\longrightarrow} \mbox{figure-eight}~\mbox{knot}~K4a1 \nonumber \\ \stackrel{(1,1)}{\longrightarrow}\mbox{{\bf Poincar\'e}}~\mbox{{\bf homology}}~\mbox{{\bf sphere}}\nonumber.\\
\end{eqnarray}

It is pictured in Fig. 3. Applying $(1,1)$-slope Dehn filling on the figure-eight knot, the sequence terminates at the spherical manifold $M_0$ that is the Poincar\'e homology sphere (also known as Poincar\'e dodecahedral space \cite{MacLachlanBook, Thurston1997}).  Conversely, drilling out along the shortest closed geodesic, as provided by SnapPy, an inverse sequence applies as follows $(h)\rightarrow (g)\rightarrow \cdots \rightarrow (b)$. At the next step, (b) transforms into a cyclic $8$-link not into the Thurston link.

Many of the links of the sequence correspond to \lq congruence' manifolds \cite{Baker2018}. A Bianchi subgroup $\Gamma_B$ of $PSL(2,\mathcal{O}_d)$ is called a \lq congruence' subgroup if there exists an ideal $I \in \mathcal{O}_d$ so that $\Gamma_B$ contains the \lq principal congruence subgroup' $\Gamma(I)=\mbox{ker}\{PSL(2,\mathcal{O}_d) \rightarrow PSL(2,\mathcal{O}_d/I)\}$. Details about these links are available on-line at \cite{Goerner}. Of course congruence Bianchi subgroups generalize congruence subgroups that played a leading role in our approach of uqc based on $\Gamma$ \cite{PlanatModular}.

The eight links under scrutiny are pictured in Fig. 3. Four of them are congruence links and are identified by a pair $(-k,I)$ meaning that they act on $\mathbb{Q}(\sqrt{-k})$ with the ideal $I$: the Thurston's link $(-3\left\langle(5+\sqrt{-3})/2  \right\rangle)$, the $6$-link L12n2256:  $(-1,\left\langle 2+0\sqrt{-1}  \right\rangle)$, the $5$-link L10n113: $(-3,\left\langle 2+0\sqrt{-3}  \right\rangle)$ and the magic link L6a5: $(-1,\left\langle(1+\sqrt{-7})/2  \right\rangle)$. The link L10n113 is also otet$10_{00027}$ and the link L12n2256 is ooct$04_{00042}$.

The $3$-manifolds that we could identify as connected to MICs are in Table 3. This becomes a cumbersome task for links larger than the $6$-link. For the link L8n7 the cardinality sequence of subgroups is $\eta_d=[63,794,23753,280162,\cdots]$, for the link L10n113 it is given by the sequence $\eta_d=[31,176,1987,7628,11682,\cdots]$ and for the $6$-link it is given by the sequence $\eta_d=[63,580,12243,94274,\cdots]$. For the latter link, SnapPy is able to build the Dirichlet domain (not shown) of symmetry group $D_6$, corresponding to the Hesse SIC while it is not the case for the one attached to the $5$-link.

We expect that the full sequence should have a meaning in the context of uqc and MICs.
 
\section{Conclusion}

It is not yet known which basic piece a practical quantum computer will be made of. The authors of \cite{PAAI2018} developed a view of universal quantum computing (uqc) based on $3$-manifolds. In a nutshell, there exists a connection between the Poincar\'e conjecture -it states that every simply connected closed $3$-manifold is homeomorphic to the $3$-sphere $S_3$- and the Bloch sphere that houses the qubits. According to this approach, the dressing of qubits in $S_3$ leads to $3$-manifolds (they have been investigated in many details by W.~P. Thurston culminating in the proof of Poincar\'e conjecture) many of them corresponding to MICs (minimal informationally complete POVMs) and the related uqc. In \cite{PlanatModular} the MICs were based on the modular group $\Gamma$ that corresponds to the trefoil knot approach in \cite{PAAI2018}. In the present paper, the $3$-manifolds under investigation derive from subgroups of Bianchi groups, a generalization of $\Gamma$. There seems to exist a nice connection between some congruence subgroups of Bianchi groups through the chain (1) and quantum computing that the author started to investigate. At a more practical level, $3$-manifolds can also be seen as (three-dimensional) quasiparticles beyond anyons \cite{Vijay2017}.

\end{document}